\documentclass[12pt]{article}
\def\date{29 August 2021}
\usepackage{amsmath, amssymb, amsfonts, amsthm}

\theoremstyle{plain}
\newtheorem{proposition}{proposition}[section]

\newtheorem{theorem}[proposition]{Theorem}
\newtheorem{claim}[proposition]{Claim}

\newtheorem{lem}[proposition]{Lemma}
\newtheorem{cor}[proposition]{Corollary}
\newtheorem{thm}[proposition]{Theorem}

\theoremstyle{definition}
\newtheorem{definition}[proposition]{Definition}

\begin{document}
\font\smallrm=cmr8

\phantom{a}\vskip .25in
\centerline{{\large \bf  FIVE-LIST-COLORING GRAPHS ON SURFACES:}}
\smallskip
\centerline{{\large\bf THE MANY FACES FAR-APART GENERALIZATION}}
\smallskip
\centerline{{\large\bf OF THOMASSEN'S THEOREM}}
\vskip.4in
\centerline{{\bf Luke Postle}%
\footnote{\texttt{lpostle@uwaterloo.ca}. Partially supported by NSERC under Discovery Grant No. 2019-04304 and the Canada Research Chairs program.}} 
\smallskip
\centerline{Department of Combinatorics and Optimization}
\centerline{University of Waterloo}
\centerline{Waterloo, ON}
\centerline{Canada N2L 3G1}
\medskip
\centerline{and}

\medskip
\centerline{{\bf Robin Thomas}%
\footnote{\texttt{thomas@math.gatech.edu}. Partially supported by NSF under
Grant No.~DMS-1202640.}}
\smallskip
\centerline{School of Mathematics}
\centerline{Georgia Institute of Technology}
\centerline{Atlanta, Georgia  30332-0160, USA}

\vskip 1in \centerline{\bf ABSTRACT}
\bigskip

{
\noindent
Let $G$ be a plane graph with $C$ the boundary of the outer face and let $(L(v):v\in V(G))$ be a family of non-empty sets.
By an \emph{$L$-coloring} of a subgraph $J$ of $G$ we mean a (proper) coloring $\phi$ of $J$  such that $\phi(v)\in L(v)$ for every vertex $v$ of $J$.
Thomassen proved that if $v_1,v_2\in V(C)$ are adjacent, $L(v_1)\ne L(v_2)$,
$|L(v)|\ge3$ for every $v\in V(C)\setminus \{v_1,v_2\}$ and $|L(v)|\ge5$ for every $v\in V(G)\setminus V(C)$,
then $G$ has an $L$-coloring.

As one final application in this last part of our series on $5$-list-coloring, we derive from all of our theory a far-reaching generalization of Thomassen's theorem, namely the generalization of Thomassen's theorem to arbitrarily many such faces provided that the faces are pairwise distance $D$ apart for some universal constant $D>0$.
}

\vfill \baselineskip 11pt \noindent \date.
\vfil\eject
\baselineskip 18pt

\section{Introduction}

\subsection{List-Coloring Graphs on Surfaces}

There exists a generalization of coloring, called \emph{list coloring}, where the vertices do not have to be colored from the same palette of colors.

We say that $L$ is a \emph{list-assignment} for a graph $G$ if $L(v)$ is a set of colors for every vertex $v$. We say $L$ is a \emph{$k$-list-assignment} if $|L(v)|\ge k$ for all $v\in V(G)$. We say that a graph $G$ has an \emph{$L$-coloring} if there exists a coloring $\phi$ such that $\phi(v)\in L(v)$ for all $v\in V(G)$. We say that a graph $G$ is \emph{$k$-choosable}, also called \emph{$k$-list-colorable}, if for every $k$-list-assignment $L$ for $G$, $G$ has an $L$-coloring. The \emph{list chromatic number} of $G$, denoted by $ch(G)$, is the minimum $k$ such that $G$ is $k$-list-colorable.

One notable difference between list coloring and ordinary coloring is that the 
Four Color Theorem~\cite{4CT1, 4CT2} does not generalize to list-coloring. 
Indeed, Voigt~\cite{Voigt} constructed a planar graph that is not $4$-choosable. 
On the other hand Thomassen~\cite{ThomPlanar} proved that every planar graph is $5$-choosable.
His proof is remarkably short and beautiful. For the sake of the inductive argument he proves the following
stronger statement.

\begin{theorem}[Thomassen]\label{Thom}
If $G$ is a plane graph with outer cycle $C$ and $P=p_1p_2$ is a path of length one in $C$ and $L$ is a list assignment for $G$ with $|L(v)|\ge 5$ for all $v\in V(G)\setminus V(C)$, $|L(v)|\ge 3$ for all $v\in V(C)\setminus V(P)$, and $|L(p_1)|=|L(p_2)|=1$ with $L(p_1)\ne L(p_2)$, then $G$ is $L$-colorable.
\end{theorem}

To state our main theorem, we need some definitions as follows.

\begin{definition}
Let $G$ be a plane graph and $f$ a face of $G$, the \emph{boundary graph} of $f$ in $G$, denoted $\partial_G f$, is the subgraph of $G$ consisting of the vertices and edges on the boundary of $f$. A \emph{chord of $f$} is an edge $e$ that is not in the boundary graph of $f$ in $G$ but whose ends are both in the boundary graph of $f$ in $G$.
\end{definition}

\begin{definition}
Let $G$ be a plane graph and $L$ a list assignment for $G$. We say a face $f$ of $G$ is \emph{restricted} under $L$ if there exists a path $P$ in $\partial_G f$ of length at most one such that $|L(v)|\ge 3$ for all $v\in V(\partial_G f)\setminus V(P)$, $|L(v)|\ge 1$ for all $v\in V(P)$ and $P$ has an $L$-coloring. In that case, we say $P$ is a \emph{special path} of $f$. 
\end{definition}

In this langauge, here is Thomassen's result restated.

\begin{theorem}[Thomassen]\label{Thom2}
Let $G$ be a plane graph and $f$ be a face of $G$. If $L$ is a list assignment for $G$ such that $f$ is restricted under $L$ and $|L(v)|\ge 5$ for all $v\in V(G)\setminus V(\partial_G f)$, then $G$ has an $L$-coloring.
\end{theorem}

Here is our main result, a far-reaching generalization of Theorem~\ref{Thom2} to arbitrarily many faces provided they are pairwise far apart.

\begin{thm}\label{MainThm}
There exists $D>0$ such that the following holds: Let $G$ be a plane graph and $\mathcal{F} = \{f_1,\ldots, f_m\}$ a set of faces of $G$ such that $d_G(\partial_G f_i, \partial_G f_j)\ge D$ for all $i\ne j\in [m]$. If $L$ is a list assignment for $G$ such that $f_i$ is restricted under $L$ for every $i\in[m]$ and $|L(v)|\ge 5$ for all $v\in V(G)\setminus \bigcup_{i\in[m]} V(\partial_G f_i)$, then $G$ has an $L$-coloring.
\end{thm}

In the previous paper in our series, we proved that Theorem~\ref{MainThm} holds when $m=2$ and $f_1,f_2$ have no special paths. Indeed that result is one of the key main results of the entire series. Hence why Theorem~\ref{MainThm} is such a vast generalization of Theorem~\ref{Thom2}. Nevertheless using the tools we have developed throughout this series, the proof of Theorem~\ref{MainThm} will be relatively short. In fact, for inductive purposes, we prove the following stronger version of Theorem~\ref{MainThm}, but first a definition.

\begin{definition}
Let $G$ be a plane graph and $L$ a list assignment for $G$. We say a subset $X$ of $V(G)$ is \emph{restricted} under $L$ if there exists a face $f$ of $G$ and a path $P$ in $\partial_G f$ of length at most one such that $|L(v)|\ge 3$ for all $v\in X$, $|L(v)|\ge 1$ for all $v\in V(P)$ and $P$ has an $L$-coloring. In that case, we say $P$ is a \emph{special path} of $X$ and $f$ is a \emph{special face} of $X$. 
\end{definition}

Here is our stronger theorem.

\begin{thm}\label{MainThm2}
There exists $D>0$ such that the following holds: Let $G$ be a plane graph and $\mathcal{X} = \{X_1,\ldots, X_m\}$ be a set of subsets of $V(G)$ such that $d(X_i.X_j)\ge D$ for all $i\ne j\in [m]$. If $L$ is a list assignment for $G$ such that $X_i$ is restricted under $L$ for every $i\in[m]$ and $|L(v)|\ge 5$ for all $v\in V(G)\setminus \bigcup_{i\in[m]} X_i$, then $G$ has an $L$-coloring.
\end{thm}

Clearly Theorem~\ref{MainThm} follows Theorem~\ref{MainThm2} by setting $X_i = V(\partial_G f_i)$ for each $i\in [m]$.

\subsection{Outline of Paper}
In Section~\ref{SecCanvases}, we recall the definition of canvas and critical canvas. In Section~\ref{SecBottleneck}, we recall the many paths bottleneck theorem. In Section~\ref{SecLinear}, we recall our bound that the number of vertices in a planar $H$-critical graph $G$ is at most linear in $|V(H)|$ and show how this implies exponential (in $r$) growth in the $r$-neighborboods of vertices of $G$ that do not intersect $H$. In Section~\ref{SecSteiner}, we recall the notion of an optimal Steiner tree and derive some of its properties. Finally, in Section~\ref{MainProof}, we prove our main technical theorem, Thoerem~\ref{MainThm2}.

\section{Canvases}\label{SecCanvases}

Let us recall the definition of canvas.

\begin{definition}[Canvas]

We say that $(G, S, L)$ is a \emph{canvas} if $G$ is a connected plane graph, $S$ is a subgraph of the boundary of the infinite face of $G$, and $L$ is a list assignment for $G$ such that $|L(v)|\ge 5$ for all $v\in V(G)\setminus V(C)$ where $C$ is the boundary of the infinite face of $G$, $|L(v)|\ge 3$ for all $v\in V(G)\setminus V(S)$, and there exists an $L$-coloring of $S$. 

If $S$ is a path that is also a subwalk of the the outer walk of $G$, then we say that $(G,S,L)$ is a \emph{path-canvas}.
\end{definition}

Hence, Thomassen's theorem restated in these terms is as follows.

\begin{thm}[Thomassen]\label{Thomassen}
If $(G,P,L)$ is a path-canvas and $|V(P)|\le 2$, then $G$ is $L$-colorable. 
\end{thm}
 
We also need the following definition of a critical canvas.

\begin{definition}[$T$-critical]

Let $G$ be a graph, $T \subseteq G$ a (not necessarily induced) subgraph of $G$ and $L$ a list assignment for $G$. For an $L$-coloring $\phi$ of $T$, we say that \emph{$\phi$ extends to an $L$-coloring of $G$} if there exists an $L$-coloring $\psi$ of $G$ such that $\phi(v)=\psi(v)$ for all $v\in V(T)$. The graph $G$ is \emph{$T$-critical with respect to the list assignment $L$} if $G \ne T$ and for every proper subgraph $G' \subset G$ such that $T \subseteq G'$, there exists an $L$-coloring of $T$ that extends to an $L$-coloring of $G'$, but does not extend to an $L$-coloring of $G$. If the list assignment is clear from the context, we shorten this and say that $G$ is $T$-critical.
\end{definition}

\begin{definition}[Critical Canvas]
We say a canvas $(G,S,L)$ is \emph{critical} if $G$ is $S$-critical with respect to the list assignment $L$.
\end{definition}

In 2007, Thomassen~\cite{ThomWheels} characterized the critical path-canvases with $|V(P)|=3$ as follows. Thomassen called such $P$-critical graphs \emph{generalized wheels}, but following the fourth paper in our series~\cite{PT4}, we call the obstructing canvases \emph{bellows}. 
  
\begin{thm}[Thomassen]\label{P3}
If $T=(G,P,L)$ is a critical path-canvas with $|V(P)|=3$, then $T$ is a bellows.
\end{thm}

A useful fact that is not too hard to work out from Thomassen's work~\cite{ThomWheels} (which we did as Corollary 5.3 in~\cite{PT4}) is the following.

\begin{lem}\label{OuterBellowsList}
If $T=(G,P,L)$ is a bellows and $v\in V(G)\setminus V(P)$ is on the outer face of $G$, then $|L(v)|=3$. 
\end{lem}

Given the structure of bellows, this implies the following.

\begin{lem}\label{BellowsListNeighbor}
If $T=(G,P,L)$ is a bellows and $v\in V(G)\setminus V(P)$, then either $|L(v)|=3$ or $v$ has a neighbor $u\in V(G)\setminus V(P)$ with $|L(u)|=3$. 
\end{lem}

We will need the following useful lemma.

\begin{lem} \label{CriticalCut}
Let $G$ be an $S$-critical graph with respect to a list assignment $L$ for $G$. Let $G_1,G_2$ be subgraphs of $G$ such that $G=G_1\cup G_2$, $G_2\cap G_1$ is a proper subgraph of $G_2$ and $S\subseteq G_1$. Then $G_2$ is $G_1\cap G_2$-critical with respect to $L$.
\end{lem}
\begin{proof}
Since $G$ is $S$-critical, every isolated vertex of $G$ belongs to $S$ and hence to $G_1$. It follows that every isolated vertex of $G_2$ belongs to $G_1\cap G_2$. Suppose for a contradiction that $G_2$ is not $G_1\cap G_2$-critical. Then, there exists an edge $e \in E(G_2) \setminus E(G_1\cap G_2)$ such
that every $L$-coloring of $G_1\cap G_2$ that extends to $G_2 \setminus e$ also extends to $G_2$. Since $S$ is a subgraph of $G_1$, we have that $e \not\in E(S)$. Since $G$ is $S$-critical, there exists an $L$-coloring $\psi$ of $S$ that
extends to an $L$-coloring $\phi$ of $G \setminus e$ but does not extend to an $L$-coloring
of $G$. Let $\psi'$ be the restriction of $\phi$ to $G_1\cap G_2$. Then $\psi'$ is an $L$-coloring of $G_1\cap G_2$ that extends to $G_2\setminus e$ but not to $G$, a contradiction.
\end{proof}

We will also need its following useful corollary.

\begin{cor} \label{CriticalSupergraph}
Let $G$ be an $S$-critical graph with respect to a list assignment $L$ for $G$. If $T$ is a proper subgraph of $G$ containing $S$, then $G$ is $T$-critical with respect to $L$.
\end{cor}
\begin{proof}
Apply Lemma~\ref{CriticalCut} with $G_1=T$ and $G_2=G$.
\end{proof}

\section{Bottleneck Theorem for Many Paths}\label{SecBottleneck}

Our first tool for the proof of Theorem~\ref{MainThm2} is the many paths bottleneck theorem, which appears in the sixth paper of our series~\cite{PT6} and is also Theorem 3.12.1 in~\cite{PostleThesis}. We state it below but first we require a few definitions.

\begin{definition}
Let $T=(G,S,L)$ be a canvas. We say a cutvertex $v$ of $G$ is \emph{essential} if whenever $v$ divides $G$ into graphs $G_1,G_2\ne G$, where $V(G_1)\cap V(G_2)=\{v\}$ and $G_1\cup G_2=G$, then $S\cap (V(G_i)\setminus \{v\})\ne \emptyset$ for all $i\in\{1,2\}$. 
\end{definition}

\begin{definition}
Let $T=(G,S,L)$ be a canvas and $C$ be the outer walk of $G$. We say a path $P'$ in $G$ is a \emph{span} if the ends of $P'$ have lists of size less than five, and there exists a path $\delta(P')$ from the ends of $P'$ in $C$ such that the interior of $\delta(P')$ contains neither a vertex of $S$ nor an essential cutvertex. We define an \emph{exterior} of $P'$, denoted by ${\rm Ext}(P')$ as the set of vertices in $\delta(P')\cup {\rm Int}(P'\cup \delta(P'))$.
\end{definition}

\begin{definition}
Let $T=(G,S,L)$ be a canvas. We say a vertex $v\in V(G)$ is \emph{superfluous} if $v\not\in V(S)$ and there exists a span $P'$ in $G$, $|V(P')|=3$ such that $v\in {\rm Ext}(P')$. 
We say a vertex is \emph{substantial} if it is not superfluous. We define the \emph{truncation} of $G$, denoted by $G^*$, to be the subgraph of $G$ induced by the substantial vertices of $G$. We define the \emph{truncated outer walk}, denoted by $C^*$, to be the outer walk of $G^*$. 
\end{definition}

\begin{lem}\label{SubNeighbors}
Let $T=(G,S,L)$ be a critical canvas. Let $R=\{v\in V(G)\setminus V(S): |L(v)|=3\}$. Then $V(G)\setminus (R\cup N(R)) \subseteq V(G^*)$.
\end{lem}
\begin{proof}
Let $v\in V(G)\setminus V(G^*)$. By definition there exists a span $P'$ of $G$ such that $v\in {\rm Ext}(P')$. Let $G_1 = G\setminus {\rm Ext}(P')$ and let $G_2 = P' \cup \delta(P') \cup {\rm Int}(P'\cup \delta(P'))$. Note it follows from the definitions of span and essential cutvertex that $S\subseteq G_1$. Also note that $G_1\cap G_2= P'$ and that $P'$ is a proper subgraph of $G_2$ since $v\in V(G_2)\setminus V(P')$. By Lemma~\ref{CriticalCut}, we have that $G_2$ is $P'$-critical with respect to $L$.

Hence $T'=(G_2, P', L)$ is a critical canvas. Since $|V(P')|=3$, we have by Theorem~\ref{P3} that $T'$ is a bellows. By Lemma~\ref{BellowsListNeighbor}, we have that $v\in N(R)\cup R$ as desired.
\end{proof}

We recall the definition of bottleneck as follows.

\begin{definition}[Bottleneck]
Let $T=(G,S,L)$ be a canvas and $C$ be the outer walk of $G$. Suppose there exists chords $U_1,U_2$ of $C$ with no end in $S$ such that $U_1$ divides $G$ into two graphs $G_1,G_1'$ and $U_2$ divides $G$ into $G_2,G_2'$ where $G_1\cap S=G_2\cap S$. Let $G'=G\setminus (G_1\setminus U_1)\setminus (G_2\setminus U_2)$. If the canvas $T'=(G',U_1\cup U_2, L)$ contains an accordion or a harmonica, call it $T''$, we say that $T''$ is a \emph{bottleneck} of $T$. 
\end{definition}

Note for the definitions of harmonica and accordion we refer the reader to~\cite{PT3} and~\cite{PT5} respectively (or to the thesis of the first author~\cite{PostleThesis} completed under the supervision of the second author.). 

Here is the many paths bottleneck theorem. Note there was a typo in the statement in~\cite{PostleThesis}, namely $G^*$ in the conclusion was inadvertently replaced with $G$ (which would be a false statement given the existence of bellows).

\begin{thm}\label{LinearManyPaths}[Linear Bottleneck Theorem: Many Paths]
If $T=(G,S,L)$ is a connected critical canvas with outer face $C$, where $S$ is the union of disjoint paths of $C$ such that there is no bottleneck $T'=(G',U_1\cup U_2,L)$ of $T$ where $d(U_1,U_2)\ge d$, then $|V(G^*)|=O(d|S|)$.
\end{thm}

Since both accordions and harmonicas with $d(U_1,U_2) \ge 3$ contain a chord of the outer face whose both ends have lists of size exactly three in $L$, we have the following corollary.

\begin{cor}\label{LinearManyPaths2}
There exists $c_1\ge 1$ such that following holds: If $T=(G,S,L)$ is a connected critical canvas with outer face $C$, where $S$ is the union of disjoint paths of $C$ such that there is no chord $uv$ of $C$ such that $|L(u)|,|L(v)| \le 3$ and $u,v\in V(G)\setminus V(S)$, then $|V(G^*)|\le c_1 \cdot |S|$.
\end{cor}

\section{A Linear Bound and Exponential Growth}\label{SecLinear}

The following theorem is a special case of Theorem 3.18 in~\cite{Hyperbolic}. We note that the proof of Theorem 3.18 in~\cite{Hyperbolic} (via Theorem 3.6 in~\cite{Hyperbolic}) uses the previous main results of this series as a black box for its proof.

\begin{thm}\label{LinearBoundPlanar}
There exists $c_2\ge 1$ such that the following holds: Let $G$ be a planar graph, $H$ a proper subgraph of $G$ and $L$ a $5$-list-assignment of $V(G)$. If $G$ is $H$-critical with respect to $L$, then $|V(G)| \le c_2\cdot |V(H)|$.
\end{thm}

Our second key tool for the proof of Theorem~\ref{MainThm2} is the following corollary that in an $H$-critical graph $G$, the $r$-neighborboods of vertices of $G$ that do not intersect $H$ have exponential growth (in $r$).

Let $G$ be a graph and $v\in V(G)$. For every integer $r\ge 0$, let $N_r(v) = \{u\in V(G): d_G(u,v) = r\}$ and $B_r(v) = \{u\in V(G): d_G(u,v) \le  r \}$. 

\begin{cor}\label{ExpGrowth}
Let $G$ be a planar graph, $H$ a proper subgraph of $G$ and $L$ a $5$-list-assignment for $G$. If $G$ is $H$-critical with respect to $L$ and $v\in V(G)\setminus V(H)$, then for all $0\le r\le d(v,H)$, $|N_r(v)|\ge 2^{\frac{r}{2\cdot c_2}}$ where $c_2$ is the constant in Theorem~\ref{LinearBoundPlanar}.
\end{cor}
\begin{proof}
We proceed by induction on $r$. If $r=0$, then $|N_r(v)|=1\ge 2^0$ as desired. 

So we assume $r > 0$. Let $G_1= G\setminus B_{r-1}(v)$ and $G_2=G[B_r(v)]$. Hence $G_1\cap G_2= G[N_r(v)]$. Since $r>0$, we have that $G_1\cap G_2$ is a proper subgraph of $G_2$. Since $r\le d(v,H)$, we find that $H \subseteq G_1$. Note that $G=G_1\cup G_2$. Hence by Lemma~\ref{CriticalCut}, we have that $G_2$ is $G_1\cap G_2$-critical with respect to $L$.

But then by Theorem~\ref{Thom}, it follows that $|V(G_1\cap G_2)| > 1$. Hence if $r \le 2\cdot c_2$, then $|N_r(v)|\ge 2 \ge 2^{\frac{r}{2c_2}}$ as desired.

So we assume $r > 2\cdot c_2$. By Theorem~\ref{LinearBoundPlanar}, we have that $|V(G_2)| \le c_2 \cdot |V(G_1)\cap V(G_2)| = c_2 \cdot |N_r(v)|$. Since $r > 2\cdot c_2$, we find that $\bigcup_{i=r-2\cdot c_2}^{r-1} N_i(v) \subseteq V(G_2)$. Hence $\sum_{i=r-2\cdot c_2}^{r-1} |N_i(v)| \le c_2 \cdot |N_r(v)|$. It follows that there exist $i$ such that $r-2\cdot c_2 \le i \le r-1$ such that $|N_i(v)| \le \frac{|N_r(v)|}{2}$. By induction, $|N_i(v)| \ge 2^{\frac{i}{2\cdot c_2}} \ge 2^{\frac{r-2\cdot c_2}{2\cdot c_2}}$. But then $|N_r(v)| \ge 2\cdot |N_i(v)| \ge 2^{\frac{r}{2\cdot c_2}}$ as desired.
\end{proof}

\section{Optimal Steiner Trees}\label{SecSteiner}

\begin{definition}
Let $G$ be a graph and $S\subseteq V(G)$. We say $T\subseteq G$ is an \emph{optimal Steiner tree} for $S$ if $T$ is a connected subgraph of $G$ such that $S\subseteq V(T)$ and subject to that $|E(T)|$ is minimized. Let $T'$ be the tree formed from $T$ by suppressing degree two vertices not in $S$. If $e\in E(T')$, we let $\pi(e)$ denote the path in $T$ between the endpoints of $e$. We say that the path $\pi(e)$ is a \emph{seam} of the tree $T$, we let $|\pi(e)|$ denote the length of $\pi(e)$, and ${\rm mid}(\pi(e))$ denote a mid-point of that path.
\end{definition}

Note we use $\pi(e)$ for the seam in the definition above, but in the rest of the paper we will just use $e$ (or $f$) for the name of the seam and hence ${\rm mid}(e)$ for its mid-point and $|e|$ for its length. We remark that the problem of finding an optimal Steiner tree is known as the Steiner tree problem in combinatorial optimization, hence the choice of name.

Our next lemma details various properties of an optimal Steiner tree.
\begin{lem}\label{Steinerization2}
If $H$ is an optimal Steiner tree of $G$ for $S$, then
\begin{enumerate}
\item[(1)] $H$ is a tree, and
\item[(2)] every leaf of $H$ is in $S$, and
\item[(3)] the number of seams of $H$ is at most $2(|S|-1)$, and
\item[(4)] for all seams $e$ of $H$, $d_G({\rm mid}(e),S) \ge \frac{|e|-1}{2}$, and
\item[(5)] for all distinct seams $e,f$ of $H$, $d_G({\rm mid}(e),{\rm mid}(f)) \ge \frac{|e|+|f|-2}{4}$.
\end{enumerate}
\end{lem}
\noindent{\bf Proof of (1).} Suppose $H$ is not a tree. Then $H$ contains a cycle $C$. Let $e\in E(C)$. Then $H-e$ is a connected subgraph of $G$ such that $S\subseteq V(H-e)$ and yet $|E(H-e)| < |E(H)|$, contradicting the minimality of $H$.
\\

\noindent{\bf Proof of (2).} Suppose not and let $v$ be a leaf of $H$ not in $S$. Then $H-v$ is a connected subgraph of $G$ such that $S\subseteq V(H-v)$ and yet $|E(H-v)| < |E(H)|$, contradicting the minimality of $H$.
\\

\noindent {\bf Proof of (3).} Let $H'$ be the tree formed from $H$ by suppressing degree two vertices not in $S$. Then $|E(H')|$ is the number of seams of $H$. Yet by (2), it follows that every leaf of $H'$ is in $S$. Moreover, no vertex in $H'$ has degree exactly two. Hence $2(|V(H')|-1) = 2|E(H')| \ge |S|-1+3(|V(H')|-|S|)$. Thus $|V(H')|\le 2|S|-1$. Hence $|E(H')|\le 2(|S|-1)$ as desired.
\\

\noindent {\bf Proof of (4).} Suppose not; that is, there exists a seam $e$ of $H$ such that $d_G({\rm mid}(e),S) < \frac{|e|-1}{2}$. Thus there exists a path $P$ from ${\rm mid}(e)$ to a vertex $v$ in $S$ with $|E(P)| < \frac{|e|-1}{2}$. Let $H_1$ be the component of $H- {\rm mid}(e)$ containing $v$ and let $u$ be the end of $e$ in $H_1$. Let $P_1$ be the path from ${\rm mid}(e)$ to $u$ in $H$. Note $|E(P_1)| \ge \frac{|e|-1}{2}$ and hence $|E(P_1)| > |E(P)|$.  Now $H'=(H\setminus P_1)\cup P$ is a connected subgraph of $G$ containing $S$ such that $|E(H')| < |E(H)|$, contradicting the minimality of $H$. 
\\

\noindent {\bf Proof of (5).} Suppose not; that is, there exist distinct seams $e, f$ of $H$ and a path $P$ from ${\rm mid}(e)$ to ${\rm mid}(f)$ with $|E(P)| < \frac{|e|+|f|-2}{4}$. We may assume without loss of generality that $|e| \ge |f|$. Let $H_1$ be the component of $H- {\rm mid}(e)$ containing ${\rm mid}(f)$ and let $u$ be the end of $e$ in $H_1$. Let $P_1$ be the path from ${\rm mid}(e)$ to $u$ in $H$. Note $|E(P_1)| \ge \frac{|e|-1}{2}$. Since $|e|\ge |f|$, we find that $|E(P_1)| \ge \frac{|e|+|f|-2}{4}$. Hence $|E(P_1)| > |E(P)|$.  Now $H'=(H\setminus P_1)\cup P$ is a connected subgraph of $G$ containing $S$ such that $|E(H')| < |E(H)|$, contradicting the minimality of $H$. \qed

\section{Main Proof}\label{MainProof}

We are now ready to prove Theorem~\ref{MainThm2}.

\begin{proof}[Proof of Theorem~\ref{MainThm2}.]
Let $c_1\ge 1$ be the constant in Corollary~\ref{LinearManyPaths2} and $c_2\ge 1$ be the constant in Theorem~\ref{LinearBoundPlanar}. Let $D\ge 720\cdot c_1\cdot c_2^2$ be an even integer large enough such that
$$2^{\frac{D-4-204\cdot c_1}{32\cdot c_1\cdot c_2}} > \frac{D-4}{4\cdot c_1}.$$ 

Suppose not. Let $G$ be a counterexample with $|V(G)|+|E(G)|$ minimized and subject to that $\sum_{v\in V(G)}|L(v)|$ minimized. Since $|V(G)|$ is minimized, it follows that $G$ is connected. For each $i\in [m]$, let $P_i$ be a special path of $X_i$ and $f_i$ a special face of $X_i$. Since $\sum_{v\in V(G)} |L(v)|$ is minimized, it follows that $P_i\ne \emptyset$ for every $i\in [m]$ and that $|L(v)|=1$ for all $v\in \bigcup_{i\in[m]} V(P_i)$. Let $P = \bigcup_{i\in [m]} P_i$. Since $|V(G)|+|E(G)|$ is minimized, we have that $G$ is $P$-critical with respect to $L$.

\begin{claim}\label{NoChord}
For every $i\in[m]$, there does not exist an edge $e\in E(G[X_i])\setminus E(\partial_G f)$.
\end{claim}
\begin{proof}
Suppose not; that is there exists $e=v_1v_2\in E(G[X_i])\setminus E(\partial_G f)$ for some $i\in [m]$. Then $e$ separates $G$ into two graphs $G_1$ and $G_2$ such that $G_1\cap G_2=e$, $G=G_1\cup G_2$ and $V(G)\setminus V(G_i)\ne\emptyset$ for each $i\in\{1,2\}$. We may assume without loss of generality that $P_i\subseteq G_1$. By the minimality of $G$, we have that $G_1$ has an $L$-coloring $\phi$. Let $L'$ be a list assignment for $G_2$ such that $L'(v_i)=\{\phi(v_i)\}$ for $i\in \{1,2\}$ and $L'(v)=L(v)$ for all $v\in V(G_2)\setminus \{v_1,v_2\}$. Let $J = \{j\in m: X_j\cap V(G_2)\ne\emptyset\}$. For each $j\in J$, let $X_j' = X_j\cap V(G_2)$ and let $f_j'$ be the face of $G_2$ containing $f_j$. For each $j\in J\setminus \{i\}$, let $P_j' = P_j$ and let $P_i' = v_1v_2$.

Now $G_2$ is a plane graph and $\mathcal{X}' = \{X_j': j\in J\}$ is a set of subsets of $V(G_2)$ such that $d_{G_2}(X_j',X_k') \ge D$ for all $j \ne k\in J$. Moreover, $L'$ is a list assignment for $G_2$ such that $X_j'$ is restricted under $L'$ for every $j\in J$ and $|L'(v)|\ge 5$ for all $v\in V(G_2)\setminus \bigcup_{j\in J} X_j'$. Since $|V(G_2)| < |V(G)|$, we have by the minimality of $G$ that $G_2$ has an $L'$-coloring $\phi'$. But then $\phi\cup \phi'$ is an $L$-coloring of $G$, a contradiction.
\end{proof}

Let $\mathcal{F}$ be the set of faces $f$ of $G$ such that $f=f_i$ for some $i\in[m]$. Let $G'$ be the graph obtained from $G$ by adding a new vertex $v_{f}$ to every face $f\in\mathcal{F}$ and adding edges from $f$ to every vertex in $V(\partial_G f)$. Let $Y=\{v_{f}: f\in \mathcal{F}\}$. Let $H$ be an optimal Steiner tree for $Y$ in $G'$, and subject to that the number of edges of $H$ incident with vertices in $Y$ is maximized. Note the latter condition implies that for every $uv\in E(P)$, $u$ and $v$ are in different components of $(H\setminus Y)\setminus E(P)$ and hence that $E(P)\cap E(H)=\emptyset$. 

Let $\mathcal{E}$ be the set of seams of $H$. Let $\mathcal{E}_1 = \{ e\in \mathcal{E}: |e| \ge 16\}$ and $\mathcal{E}_2 = \mathcal{E}\setminus \mathcal{E}_1$. Let $R=\{v\in V(G)\setminus (V(S)\cup V(H)): |L(v)|=3\}$ and $Z=V(G)\setminus (R\cup N(R))$. 

\begin{claim}\label{ZLower}
$|Z| \ge \frac{D-4}{2} \cdot m$.
\end{claim}
\begin{proof}
For each $i\in [m]$, let $Q_i=v_{i,0}v_{i,1}\ldots v_{i, k_i}$ be a shortest path in $G$ from $X_i$ to $\bigcup_{j\ne i} X_j$. Since $d_G(X_i,X_j)\ge D$ for all $i\ne j\in [m]$, we find that $k_i \ge D$ for all $i\in [m]$. For each $i\in [m]$, let $Q_i'=v_{i,2}v_{i,3}\ldots v_{i, \frac{D-2}{2}}$; note $|Q_i'| = \frac{D-4}{2}$ and moreover $V(Q_i') \subseteq Z$. Since $d_G(X_i,X_j)\ge D$ for all $i\ne j\in [m]$, we find that $V(Q_i')\cap V(Q_j') = \emptyset$ for all $i\ne j \in [m]$. Thus $|Z| \ge \sum_{i\in [m]}|V(Q_i')| \ge \frac{D-4}{2} \cdot m$ as desired.
\end{proof}

Let $G_0$ be the graph obtained from $G$ by cutting along the seams of $H$ as follows. Let $H'$ be the graph obtained from $H$ by suppressing degree two vertices. Every seam $e\in \mathcal{E}$ becomes two copies $e^L$ and $e^R$ on the left and right sides respectively of the cut; if $v\in V(H')$ is incident with seams $e_1, \ldots, e_k$ where $k={\rm deg}_H(v)$, then $v$ is split into $k$ copies $v_1,\ldots,v_k$ where $v_i$ is the end of the seams $e_{i}^R$ and $e_{(i+1) \mod k}^L$. Note that all of $Y$ is now on the boundary of a new face. We then delete the vertices of $Y$. Hence all vertices in $(H\setminus Y) \cup \bigcup_{i\in [m]} X_i$ (or rather their copies) are incident with the same face of $G_0$, call it $f_0$. We may assume without loss of generality that $f_0$ is the infinite face of $G_0$. Note that $G_0$ is connected since $G$ is connected.

We define a mapping $\rho$ from $V(G_0)$ to $V(G)$ as follows: If $v\in V(H')$, we let $\rho(v_i) = v$ for all $i\in [{\rm deg}_H(v_i)]$. If $v\in V(H)\setminus V(H')$, we let $\rho(v^L) = \rho(v^R) = v$. Finally if $v\in V(G)\setminus V(H)$, we let $\rho(v) = v$. We then extend this mapping to edges $e=uv$ of $G_0$ by letting $\rho(e)=\rho(u)\rho(v)$ and subgraphs $G_0'$ of $G_0$ by letting $\rho(G_0') = \{ \rho(v): v\in V(G_0')\} \cup \{\rho(uv): uv\in E(G_0')\}$.

We define a list assignment $L_0$ of $G_0$ by setting $L_0(v) = L(\rho(v))$ for all $v\in V(G_0)$. Let $S=P\cup (H\setminus Y)$. Note that $S$ is a subgraph of $G$. We let 
$$S_0 = \{v\in V(G_0): \rho(v)\in V(S)\} \cup \{e=uv\in E(G_0): \rho(u)\rho(v)\in E(S)\}.$$ 
Note that $S_0$ is a subgraph of the boundary of the infinite face of $G_0$. 

Note that every $L$-coloring $\phi$ of $S$ corresponds to an $L_0$-coloring $\phi_0$ of $S_0$ by letting $\phi_0(v) = \phi(\rho(v))$ for all $v\in V(S_0)$. 

\begin{claim}\label{S0Colorable}
$S$ has an $L$-coloring.
\end{claim}
\begin{proof}
Let $\phi$ be an $L$-coloring of $P$. Note such a coloring exists by assumption. Since $|L(v)|=1$ for all $v\in V(P)$, this coloring is unique. Recall that for every $uv\in E(P)$, $u$ and $v$ are in different components of $(H\setminus Y)\setminus E(P)$. Let $W$ be a component of $S\setminus V(P)$. Let  $L'(v)= L(v) \setminus \{\phi(u):u\in N_S(v)\cap V(P)\}$. Since $D\ge 3$, we have that $|N_S(v)\cap V(P)| \le 1$ for all $v\in S\setminus V(P)$. Thus $|L'(v)| \ge |L(v)|-1 \ge 2$ for all $v\in W$. Moreover, $W$ is a tree. Since $|L'(v)|\ge 2$ for all $v\in V(W)$, there exists an $L'$-coloring $\phi_W$ of $W$. But then $\phi\cup \bigcup_W \phi_W$ is an $L$-coloring of $S$ as desired.
\end{proof}

Since $S$ has an $L$-coloring by Claim~\ref{S0Colorable}, it follows that there exists an $L_0$-coloring of $S_0$. Thus $T_0=(G_0,S_0,L_0)$ is a canvas. Furthermore, since $S$ has an $L$-coloring of $G$, it follows that $S$ is a proper subgraph of $G$ as there does not exist an $L$-coloring of $G$. Since $G$ is $P$-critical with respect to $L$, this implies by Corollary~\ref{CriticalSupergraph} that $G$ is $S$-critical with respect to $L$. In addition, it also implies that $S_0$ is a proper subgraph of $G_0$. 

\begin{claim}\label{T0Critical}
$T_0$ is critical.
\end{claim}
\begin{proof}
Suppose not. Since $G_0$ is connected, it follows that there exists an edge $e=uv\in E(G_0)\setminus E(S_0)$ such that every $L_0$-coloring of $S_0$ that extends to $G_0 \setminus e$ also extends to $G_0$. Sinc $uv\in E(G_0) \setminus E(S_0)$, it follows that if $u'v'\in E(G_0)$ such that $\rho(u')\rho(v')\in E(G)$, then $\{u,v\}=\{u',v'\}$. Hence $\rho(G_0\setminus e) = G\setminus e'$. Let $e'=\rho(u)\rho(v)$. Thus $e' \in E(G)\setminus E(S)$. Since $G$ is $S$-critical, we find that there exists an $L$-coloring $\psi$ of $S$ such that $\psi$ extends to an $L$-coloring $\phi$ of $G\setminus e$ but not to an $L$-coloring of $G$. Let $\psi_0(v) = \psi(\rho(v))$ for all $v\in V(S_0)$ and let $\phi_0(v) = \phi(\rho(v))$ for all $v\in V(G_0\setminus e)$. Hence $\psi_0$ extends to $\phi_0$ and hence to an $L_0$-coloring of $G_0$. But then $\phi_0(u)\ne \phi_0(v)$. Hence $\phi(\rho(u))\ne \phi(\rho(v))$ and so $\phi$ is also an $L$-coloring of $G$, a contradiction.
\end{proof}

Let $R_0 = \{v\in V(G_0)\setminus V(S_0): |L_0(v)|=3\}$. Note that $R_0 = \{v\in V(G_0)\setminus V(S_0): \rho(v)\in R \}$. Let $Z_0=V(G_0)\setminus (R_0\cup N_{G_0}(R_0))$. Note that $Z_0 = \{v\in V(G_0): \rho(v)\in Z\}$. Moreover, $|R| = |R_0|$ and $|Z| \le |Z_0|$.

\begin{claim}\label{T0NoChord}
There is no chord $uv$ of $f_0$ such that $|L_0(u)|,|L_0(v)| \le 3$ and $u,v\in V(G_0)\setminus V(S_0)$.
\end{claim}
\begin{proof}
Suppose not. Let $uv$ be a chord of $f_0$ such that $|L_0(u)|,|L_0(v)| \le 3$ and $u,v\in V(G_0)\setminus V(S_0)$. Thus $u,v \in R_0$. Hence $\rho(u),\rho(v)\in R$. Thus there exists $i,j\in [m]$ with $\rho(u)\in X_i$ and $\rho(v)\in X_j$. Since $uv\in E(G_0)$, we have that $e:=\rho(u)\rho(v)\in E(G)$. Thus $d_G(X_i,X_j) \le 1$. Since $D \ge 2$, it follows that $i=j$. Since $uv$ is a chord of $f_0$, it follows that $e$ is a chord of $f_i$. Hence by the definition of chord, we have that $e\in E(G[X_i])\setminus E(\partial_G f)$, contradicting Claim~\ref{NoChord}. 
\end{proof}

By Claims~\ref{T0Critical} and~\ref{T0NoChord}, and since $S_0$ is the disjoint union of paths of $f_0$, we have by Corollary~\ref{LinearManyPaths2} that $|V(G_0^*)| \le c_1 \cdot |S_0|$. Since $T_0$ is critical by Claim~\ref{T0Critical}, we have by Lemma~\ref{SubNeighbors} that $Z_0 \subseteq V(G_0^*)$. Hence by Claim~\ref{ZLower}, we have that
$$\frac{D-4}{2} \cdot m \le |Z| \le |Z_0| \le |V(G_0^*)| \le c_1 \cdot |S_0|.$$
Yet by Lemma~\ref{Steinerization2}(3), we have that $|\mathcal{E}|\le 2m$ and hence
$$|S_0| \le |V(P)|+2\cdot \sum_{e\in \mathcal{E}} (|e|+1) \le 2m + 34\cdot |\mathcal{E}|+\sum_{e\in\mathcal{E}_1} |e|\le 70m+\sum_{e\in\mathcal{E}_1} |e|.$$
Thus
$$\sum_{e\in\mathcal{E}_1} |e| \ge \left(\frac{D-4}{2\cdot c_1}-70\right) \cdot m.$$
Let
$$x:= \frac{\sum_{e\in\mathcal{E}_1} |e|}{|\mathcal{E}_1|}.$$
Let $D'= \frac{D-4}{4\cdot c_1}-35$. We note that since $D\ge 720\cdot c_1\cdot c_2^2$ and $c_1,c_2 \ge 1$, we have that $D' \ge 144 \cdot c_2^2$. Since $|\mathcal{E}_1| \le |\mathcal{E}|\le 2m$, we find that
$$x\ge \frac{\sum_{e\in\mathcal{E}_1} |e|}{|\mathcal{E}|}\ge \frac{\sum_{e\in\mathcal{E}_1} |e|}{2m} \ge \frac{D-4}{4\cdot c_1}-35=D'.$$
Moreover,
$$|Z| \le 70m+\sum_{e\in\mathcal{E}_1} |e| = 70m+x\cdot |\mathcal{E}_1| \le (70+2x)m.$$

For each $e\in \mathcal{E}_1$, let $r_e = \left\lceil \frac{|e|}{4}-4 \right\rceil$ and let $A_e = N_{r_e}({\rm mid}(e))$. Note that $r_e \le \frac{|e|-1}{4}-3$. Recall that for each $e\in \mathcal{E}_1$, we have by definition of $\mathcal{E}_1$ that $|e|\ge 16$ and hence $r_e\ge 0$. By Lemma~\ref{Steinerization2}(4) for each $e\in \mathcal{E}_1$, we have that $d_{G'}({\rm mid}(e),Y) \ge \frac{|e|-1}{4}$. Yet for every $r\in R$, $d_{G'}(r,Y)=1$. Thus $d_G({\rm mid}(e),R\cup N(R)) \ge \frac{|e|-1}{4} - 2 > \lceil r_e \rceil$. Thus for each $e\in \mathcal{E}_1$, we have that 
$$A_e \subseteq V(G)\setminus (R\cup N(R)) = Z.$$ 

Furthermore, by Lemma~\ref{Steinerization2}(4), we have that $A_e \cap A_f = \emptyset$ for all distinct seams $e\ne f\in \mathcal{E}_1$. By Corollary~\ref{ExpGrowth}, we have that for each $e\in \mathcal{E}_1$
$$|A_e| = |N_{r_e}({\rm mid}(e))| \ge 2^{\frac{r_e}{2\cdot c_2}} \ge 2^{\frac{|e|-16}{8\cdot c_2}}.$$ 
Thus
$$|Z| \ge \sum_{e\in \mathcal{E}_1} |A_e| \ge \sum_{e\in \mathcal{E}_1} 2^{\frac{|e|-16}{8\cdot c_2}}.$$
By the convexity of the exponential function, we find that
$$|Z| \ge |\mathcal{E}_1| \cdot 2^{\frac{x-16}{8\cdot c_2}}.$$
Combining the upper and lower bounds on $|Z|$, we find that
$$ (70+2x)m \le |Z| \le |\mathcal{E}_1| \cdot 2^{\frac{x-16}{8\cdot c_2}} \le 2m \cdot 2^{\frac{x-16}{8\cdot c_2}},$$
and hence
$$35+x \le 2^{\frac{x-16}{8\cdot c_2}}.$$
Let $f(y) = 2^{\frac{y-16}{8\cdot c_2}} - (35+y)$. Thus $f(x)\le 0$. Note that $f'(y) = \frac{\ln 2}{8\cdot c_2} \cdot 2^{\frac{y-16}{8\cdot c_2}} - 1$. For every $y\ge D'$, we have that $y\ge D' \ge 144\cdot c_2^2$ and hence $2^{\frac{y-16}{8\cdot c_2}} \ge 2^{\frac{D'-16}{8\cdot c_2}} \ge 2^{16\cdot c_2} \ge 16\cdot c_2$; thus $f'(y) \ge f'(D') \ge 0$. So $f$ is increasing for every $y\ge D'$. Since $x\ge D'$, we find that 
$$f(x) \ge f(D') = f\left(\frac{D-4}{4\cdot c_1}-35\right) = 2^{\frac{D-4-204\cdot c_1}{32\cdot c_1\cdot c_2}} - \frac{D-4}{4\cdot c_1} > 0,$$
where the last inequality follows by the assumption on $D$. Yet from before, we have that $f(x) \le 0$, a contradiction.
\end{proof}

\subsubsection*{Acknowledgements.} 
The work for this research was done in 2013-2014 while the first author was a postdoc at Emory University. At the time, we opted not to release it as it relies on the other eight papers of this series, which are still getting published to this today (see~\cite{PT1, PT2, PT3, PT4}). The contents of those eight papers as well as our general paper~\cite{Hyperbolic} can essentially be found in the thesis of the first author~\cite{PostleThesis}, completed under the supervision of the second author. This work however is not contained there and so in the interests of the community, the sole surviving author has decided to release it. The first author did announce the result in a talk at the University of Waterloo in 2016 (https://uwaterloo.ca/combinatorics-and-optimization/events/graph-theory-luke-postle) with a proof as described in this paper.

The first author thanks Joan Hutchinson for insightful discussions on the topic during CanaDAM 2013.

%

\begin{thebibliography}{99}

\def\JCTB{{\it J.~Combin.\ Theory Ser.\ B}}
\def\CMUC{{\it Comment. Math. Univ. Carol.}}
\def\TAMS{{\it Trans.\ Amer.\ Math.\ Soc.}}
\def\JAMS{{\it J.~Amer.\ Math.\ Soc.}}
\def\PAMS{{\it Proc. Amer. Math. Soc.}}
\def\DM{{\it Discrete Math.}}
\def\CM{{\it Contemporary Math.}}
\def\GC{{\it Graphs and Combin.}}
\def\COM{{\it Combinatorica}}
\def\JGT{{\it J.~Graph Theory}}
\def\JAlgorithms{{\it J.~Algorithms}}
\def\SIAMDM{{\it SIAM J.~Disc.\ Math.}}
\def\CPC{{\it Combinatorics, Probability and Computing}}
\def\EJC{Electron.\ J.~Combin.}

\bibitem{4CT1} K. Appel and W. Haken, Every planar map is four colorable, Part I: discharging, Illinois J. of Math. 21 (1977), 429--490.
\bibitem{4CT2} K. Appel, W. Haken, J. Koch, Every planar map is four colorable, Part II: reducibility, Illinois J. of Math. 21 (1977), 491--567.
\bibitem{PostleThesis} L. Postle, 5-list-coloring graphs on surfaces, Ph.D. Dissertation, Georgia Institute of Technology, 2012.
\bibitem{PT1} L. Postle and R. Thomas, Five-List-Coloring Graphs on Surfaces I. Two Lists of Size Two in Planar Graphs. Journal of Combinatorial Theory Ser. B 111 (2015), pp. 234--241.
\bibitem{PT2} L. Postle and R. Thomas, Five-List-Coloring Graphs on Surfaces II. A Linear Bound for Critical Graphs in a Disk. Journal of Combinatorial Theory Ser. B 119 (2016), pp. 42--65.
\bibitem{PT3} L. Postle, R. Thomas. Five-List-Coloring Graphs on Surfaces III. One List of Size One and One List of Size Two, Journal of Combinatorial Theory Ser. B 128 (2018), pp. 1--16.
\bibitem{PT4} L. Postle, R. Thomas. Five-List-Coloring Graphs on Surfaces IV. Two Lists of Size One, manuscript.
\bibitem{PT5} L. Postle, R. Thomas. Five-List-Coloring Graphs on Surfaces V. Two Precolored Edges, manuscript.
\bibitem{PT6} L. Postle, R. Thomas. Five-List-Coloring Graphs on Surfaces VI. A Linear Bound for Paths, manuscript.
\bibitem{Hyperbolic} L. Postle, R. Thomas, Hyperbolic families and coloring graphs on surfaces, Trans. Amer. Math. Soc. Ser. B 5 (2018), 167--221.
\bibitem{ThomPlanar} C. Thomassen, Every planar graph is $5$-choosable, J. Combin. Theory Ser. B 62 (1994), 180--181.
\bibitem{ThomWheels} C. Thomassen, Exponentially many 5-list-colorings of planar graphs, J. Combin. Theory Ser. B 97 (2007), 571--583.
\bibitem{Voigt} M. Voigt, List colourings of planar graphs, Discrete Mathematics 120 (1993) 215--219.

\end{thebibliography}

\baselineskip 11pt
\vfill
\noindent
This material is based upon work supported by the National Science Foundation.
Any opinions, findings, and conclusions or
recommendations expressed in this material are those of the authors and do
not necessarily reflect the views of the National Science Foundation.
\eject

\end{document}